# ESTIMATING A CONCAVE DISTRIBUTION FUNCTION FROM DATA CORRUPTED WITH ADDITIVE NOISE

BY GEURT JONGBLOED AND FRANK H. VAN DER MEULEN

*Delft University of Technology*

We consider two nonparametric procedures for estimating a concave distribution function based on data corrupted with additive noise generated by a bounded decreasing density on $(0, \infty)$. For the maximum likelihood (ML) estimator and least squares (LS) estimator, we state qualitative properties, prove consistency and propose a computational algorithm. For the LS estimator and its derivative, we also derive the pointwise asymptotic distribution. Moreover, the rate $n^{-2/5}$ achieved by the LS estimator is shown to be minimax for estimating the distribution function at a fixed point.

**1. Introduction.** Let $X_1, X_2, \ldots$ be an i.i.d. sequence of random variables with unknown distribution function $F$. Moreover, let $\varepsilon_1, \varepsilon_2, \ldots$ be an i.i.d. sequence of random variables, independent of the $X_i$'s, with known probability density function $k$. We want to estimate the distribution function $F$, based on data $Z_1, Z_2, \ldots, Z_n$, where $Z_i = X_i + \varepsilon_i$. In other words, we wish to estimate $F$ based on a sample from the density

$$g_F(z) = \int_{\mathbb{R}} k(z-x)\,dF(x). \tag{1}$$

Since $g_F$ is the convolution of the unknown distribution function with the (known) density $k$, the problem of estimating aspects of the distribution function $F$ based on a sample from $g_F$ is known as a deconvolution problem.

Deconvolution problems were studied quite extensively during the past two decades. Given a class $\mathcal{F}$ of distribution functions $F$, one can qualitatively state that the smoother the noise density $k$, the worse the optimal estimation rate for $F$. See Fan (1991). Alternatively, given a noise density $k$, it is obvious that the smaller the class of distribution functions $\mathcal{F}$, the better the optimal estimation rate for $F$.

One popular approach to this estimation problem is based on kernel smoothing and Fourier methods [see, e.g., Carroll and Hall (1988) and Delaigle and Hall (2006)]. These estimators can achieve optimal rates of convergence under a wide range of smoothness assumptions. A characteristic feature of this approach is the need for a bandwidth, preferably chosen in an asymptotically optimal way. Many methods have been developed to determine such a bandwidth [see, e.g., Stefanski and Carroll (1990) and Delaigle and Gijbels (2004)]. Another popular approach is based on wavelets [see, e.g., Pensky and Vidakovic (1999)]. For both Fourier inversion methods and wavelet methods it is difficult to incorporate shape constraints on the distribution of interest in the estimation procedure. For example, density estimates can easily become negative.

Another method that can be employed to estimate the distribution function $F$ is maximum likelihood. Based on the density (1) of $Z_i$, the log likelihood of a density $g$ (or equivalent distribution function $F$) is easily computed. A maximum likelihood estimator is then defined as the maximizer of the log likelihood function over an appropriate class of distribution functions. See, for example, Groeneboom and Wellner (1992) for the case where it is maximized over the class of all distribution functions on $[0, \infty)$. Another general method to estimate $F$ is least squares. Based on a naive estimator of $F$ outside the class $\mathcal{F}$ of distribution functions of interest, this estimator is defined as the minimizer of the $L_2$ distance to this naive estimator over the class of interest. Typically, maximum likelihood and least squares estimators do not require a bandwidth. Moreover, shape constraints can quite naturally be imposed on the estimator by restricting the feasible set of distribution functions in their definition. This in contrast to the aforementioned kernel and wavelet based methods of estimation.

In this paper we estimate the distribution function $F$ under the assumption that it is concave. More precisely, we assume $F$ to belong to the class

(2) $\quad \mathcal{F} := \{F | F \text{ is a concave distribution function on } [0, \infty)\}$.

We restrict the convolution kernel $k$ to the class of convolution kernels

(3)
$$\mathcal{K} = \{k : [0, \infty) \to [0, \infty) : k \text{ is a}$$

bounded and decreasing probability density$\}$.

However, as pointed out in side remarks, the existence, characterization and consistency results for the maximum likelihood estimator can be extended to more general classes of kernel functions at the cost of extra technicalities.

Our initial motivation to study nonparametric estimators for shape-constrained distribution functions in deconvolution models was the financial application studied in Jongbloed, van der Meulen and van der Vaart (2005).

There, we find the problem of recovering a unimodal distribution from data corrupted with additive noise with a smooth density. The current setting with decreasing kernel $k$ is too restrictive to be applicable in that context. However, in this simplified model we can obtain asymptotic results for the LS estimator. These are of independent interest. To our knowledge, this paper is the second setting where the so-called Groeneboom distribution described in Groeneboom, Jongbloed and Wellner (2001a) appears in the limit. The first setting is that of estimating a convex decreasing density studied in Groeneboom, Jongbloed and Wellner (2001b). In both situations, the rescaling rate of the estimator is $n^{2/5}$. We expect that the role played by Chernoff's distribution [Chernoff (1964)] in situations with cube root $n$ asymptotics [Kim and Pollard (1990)] is played by the Groeneboom distribution in situations with $n^{2/5}$ asymptotics. Examples of other estimation problems where we expect this to happen are that of estimating a log concave density [Dümbgen and Rufibach (2004)] and that of estimating a concave distribution function from current status data. (We conjecture that the maximum likelihood estimator has the same asymptotics as the least squares estimator in the setting of this paper.)

In Section 2 we define two nonparametric estimators for the concave distribution function $F$: the maximum likelihood estimator and a least squares estimator. The consistency of both estimators is proved in Section 3. Computational issues of the estimators are addressed in Section 4. Subsequently, we derive an asymptotic local minimax lower bound on the optimal estimation rate for $F(x_0)$ and $f(x_0)$ in Section 5. In Section 6 we derive the asymptotic distribution of the random vector $(\tilde{F}_n(x_0), \tilde{f}_n(x_0))$. It turns out that the asymptotic variance of the LS estimator depends on the functions $k$ and $f$ in exactly the same way as the minimax lower bound of Section 5.

**2. Two nonparametric estimators: definition and properties.** In this section we define two nonparametric estimators for $F$: the maximum likelihood (ML) and least squares (LS) estimators. In the context of convex density estimation, Groeneboom, Jongbloed and Wellner (2001b) show that the ML and LS estimators have the same asymptotic pointwise behavior. The least squares estimator, however, is much more tractable to study both from an algorithmic and asymptotic point of view. The same phenomenon will be seen to occur in the deconvolution setting of this paper.

2.1. *Maximum likelihood.* Let

$$\mathcal{Z}_n = \{Z_1, \ldots, Z_n\} \tag{4}$$

be the set of observations. Denoting by $\mathbb{G}_n$ the empirical distribution function of $\mathcal{Z}_n$, the log-likelihood function evaluated at a distribution function

$F$ is given by

$$l_n(F) = \int_{\mathbb{R}} \log g_F(z) \, d\mathbb{G}_n(z), \tag{5}$$

where $g_F$ is defined as the convolution of $k$ and $F$: $g_F(z) = \int_{[0,\infty)} k(z - x) \, dF(x)$. In Groeneboom and Wellner (1992) it is shown that the maximizer of this function over the class of *all* distribution functions is a discrete distribution function with mass concentrated at the observed data points. We show that the maximum likelihood estimator of a concave distribution function $F$, based on a sample of size $n$ from $g_F$, is a proper piecewise linear distribution function that can only have changes of slope at the observed data points. We also establish a characterization of the estimator in terms of inequalities.

Define the set $\mathcal{F}_{basis} := \{F_\theta \mid \theta > 0\}$ by

$$F_\theta(x) = \frac{x}{\theta} \mathbf{1}_{[0,\theta]}(x) + \mathbf{1}_{(\theta,\infty)}, \qquad \theta > 0 \ (x \in \mathbb{R}), \tag{6}$$

that is, $F_\theta$ is the distribution function of a uniformly distributed random variable on $[0,\theta]$. Any $F \in \mathcal{F}$ can be written as a mixture of elements from $\mathcal{F}_{basis}$: there exists a probability measure $\mu = \mu_F$ on $[0,\infty)$ such that $F = \int_{[0,\infty)} F_\theta \, d\mu_F(\theta)$. In fact, $d\mu_F(\theta) = -\theta \, dF'(\theta)$. This implies

$$g_F(x) = \int_{[0,\infty)} \int_{[0,\infty)} k(x - u) \, dF_\theta(u) \, d\mu_F(\theta) = \int_{[0,\infty)} g_\theta(x) \, d\mu_F(\theta),$$

where

$$g_\theta(x) := \int_{[0,\infty)} k(x - u) \, dF_\theta(u)$$
$$= \frac{1}{\theta}(K(x) - K(x - \theta)), \qquad \theta > 0 \ (x \in \mathbb{R}). \tag{7}$$

($K$ denotes the primitive of $k$.) Thus we can reformulate the maximum-likelihood problem as to maximize $l_n(g) = \int \log g(x) \, d\mathbb{G}_n(x)$ over $\mathcal{G}$, where

$$\mathcal{G} := \left\{ g \mid g(\cdot) = \int_{[0,\infty)} g_\theta(\cdot) \, d\mu(\theta) \text{ for some probability measure } \mu \text{ on } [0,\infty) \right\}.$$

Once we know the mixing probability measure $\hat{\mu}_n$ corresponding to the maximizer $\hat{g}_n$, the maximum-likelihood estimator for $F$ is given by $\hat{F}_n = \int F_\theta \, d\hat{\mu}_n(\theta)$.

THEOREM 2.1. *Let $k \in \mathcal{K}$ as defined in* (3). *Then a maximizer $\hat{F}_n$ of* (5) *over the class of all concave distribution functions on $[0,\infty)$ exists and can be chosen to be a piecewise linear distribution function with bend points concentrated on the set of observations $\mathcal{Z}_n$.*

PROOF. We start by showing that if $\hat{F}_n$ exists, there is a version that is piecewise linear with bend points concentrated on $\{Z_1, \ldots, Z_n\}$. Consider an arbitrary concave distribution function $F$ and its linearly interpolated version (between the observed $Z_i$'s) $\bar{F}$. Then, writing $Z_{(0)} = 0$, we get for each $i$

$$g_F(Z_{(i)}) = \sum_{j=1}^{i} \int_{Z_{(j-1)}}^{Z_{(j)}} k(Z_{(i)} - y) \, dF(y)$$

(8)

$$\leq \sum_{j=1}^{i} \int_{Z_{(j-1)}}^{Z_{(j)}} k(Z_{(i)} - y) \, d\bar{F}(y) = g_{\bar{F}}(Z_{(i)})$$

implying that $l_n(F) \leq l_n(\bar{F})$. Inequality (8) holds because we can write for each summand (treating the $Z_{(i)}$'s as fixed and denoting the distribution of a uniformly distributed random variable $U$ on $[0,1]$ by $J$)

$$\int_{Z_{(j-1)}}^{Z_{(j)}} k(Z_{(i)} - y) \, dF(y) = E_F k(Z_{(i)} - Y) 1_{(Z_{(j-1)}, Z_{(j)}]}(Y)$$

$$= E_J k(Z_{(i)} - F^{-1}(U)) 1_{(Z_{(j-1)}, Z_{(j)}]}(F^{-1}(U))$$

$$\leq E_J k(Z_{(i)} - \bar{F}^{-1}(U)) 1_{(Z_{(j-1)}, Z_{(j)}]}(\bar{F}^{-1}(U))$$

$$= E_{\bar{F}} k(Z_{(i)} - Y) 1_{(Z_{(j-1)}, Z_{(j)}]}(Y)$$

$$= \int_{Z_{(j-1)}}^{Z_{(j)}} k(Z_{(i)} - y) \, d\bar{F}(y).$$

Here we use that $F^{-1}(u) \in (Z_{(j-1)}, Z_{(j)}] \iff \bar{F}^{-1}(u) \in (Z_{(j-1)}, Z_{(j)}]$ and that for each $u \in (0,1)$, $F^{-1}(u) \leq \bar{F}^{-1}(u)$ implying that $k(Z_{(i)} - F^{-1}(u)) \leq k(Z_{(i)} - \bar{F}^{-1}(u))$.

To show existence of $\hat{F}_n$, we only have to consider distribution functions having bend points at the observations and these can be parameterized as follows:

$$F = \sum_{j=1}^{n} \tau_j F_{Z_j}$$

with $\tau \in \Xi = \left\{ \tau \in \mathbb{R}^n : 0 \leq \tau_j \leq 1 \text{ for } 1 \leq j \leq n \text{ and } \sum_{j=1}^{n} \tau_j = 1 \right\}$.

Expressed in terms of $\tau$, the log likelihood function becomes $n^{-1} \times \sum_{i=1}^{n} \log(\sum_{j=1}^{n} \tau_j g_{Z_j}(Z_i))$, which is a concave function that attains a finite value for some feasible $\tau$. Since $\Xi$ is compact, existence follows. $\square$

REMARK 2.2. Existence and piecewise linearity with at most $n$ changes of slope of $\hat{F}_n$ can also be proved under the less-restrictive assumption that $k$ should be initially nondecreasing on $\mathbb{R}$, that is under the assumption that there exists a constant $M \in \mathbb{R}$ such that $k$ is nondecreasing on $(-\infty, M)$. In that situation we should allow $\hat{F}_n$ to have a point mass at zero. This implies that $\mathcal{F}_{basis}$ should be augmented with the function $\mathbf{1}_{[0,\infty)}$. In this more general setting, the bend points of the MLE can be outside the set of observed data points.

THEOREM 2.3 (Characterization of the MLE). *The* (*piecewise linear*) *distribution function $F$ maximizes* (5) *over the class $\mathcal{F}$ if and only if*

$$\int \frac{g_\theta(z)}{g_F(z)} d\mathbb{G}_n(z) \le 1 \tag{9}$$

*for all $\theta > 0$. Here $g_\theta$ is as defined in* (7). *In fact, equality holds for those $\theta$ that belong to the set of bend points of $F$.*

PROOF. First necessity. Suppose $F$ maximizes the log likelihood. Then, for all $\theta > 0$ and $\varepsilon \in [0, 1]$,

$$F + \varepsilon(F_\theta - F) \in \mathcal{F} \Rightarrow \lim_{\varepsilon \downarrow 0} \varepsilon^{-1}(l_n(F + \varepsilon(F_\theta - F)) - l_n(F)) \le 0. \tag{10}$$

Writing out this limit gives (9). That the inequality actually is an equality for those points where $\mu_F(\{\theta\}) > 0$ follows immediately upon noting that for those points $F + \varepsilon(F_\theta - F) \in \mathcal{F}$ also for small negative values of $\varepsilon$.

For sufficiency, let $\tilde{F} = \int F_\theta \, d\tilde{\mu}(\theta)$ be an arbitrary (sub-)distribution function in $\mathcal{F}$. Then,

$$\begin{aligned}
l_n(\tilde{F}) - l_n(F) &= \int \log \frac{\tilde{g}(z)}{g_F(z)} d\mathbb{G}_n(z) \le \int \left(\frac{\tilde{g}(z)}{g_F(z)} - 1\right) d\mathbb{G}_n(z) \\
&= \int \frac{1}{g_F(z)} \int g_\theta(z) \, d\tilde{\mu}(\theta) \, d\mathbb{G}_n(z) - 1 \\
&= \int \left(\int \frac{g_\theta(z)}{g_F(z)} d\mathbb{G}_n(z)\right) d\tilde{\mu}(\theta) - 1 \le 0. \qquad \square
\end{aligned}$$

2.2. *Least squares.* We now turn to an alternative nonparametric estimator for $F$, the least squares (LS) estimator. In order to define this estimator we need a "type of inverse" for the kernel $k$. In Lemma 2.4 we will prove that under mild conditions there exists a function $p$, such that $p * k(x) = id_+(x) := x\mathbf{1}_{[0,\infty)}(x)$. We now explain how we can use this result to define a least squares estimator. First note that

$$p * g(x) = (p * k) * dF(x) = (id_+ * dF)(x) = \int_0^x F(u) \, du,$$

which implies that the survival function of the random variable $X$, defined by $s = 1 - F$, satisfies

$$s(x) := U'(x) \qquad \text{with } U(x) := x - (p * g)(x).$$

Define an empirical estimate of $U$ by

$$U_n(x) = x - (p * dG_n)(x),$$

and denote the class of survival functions associated with $\mathcal{F}$ by

$$\mathcal{S} = \{s \in L^2[0, \infty) : s \text{ is nonnegative, convex, decreasing and } s(0) \in (0, 1]\}.$$

We would like to define the LS estimator $\hat{s}_n$ by $\arg\min_{s \in \mathcal{S}} Q_n(s)$, where

(11) $$Q_n(s) = \frac{1}{2} \int_0^\infty s(x)^2 \, dx - \int_0^\infty s(x) \, dU_n(x).$$

This definition is motivated by considering the $L^2$-distance between $s$ and (the nonexistent) $U_n'$. In the decomposition

$$\int (s(x) - U_n'(x))^2 \, dx = \int s(x)^2 \, dx - 2 \int s(x) U_n'(x) \, dx + \int U_n'(x)^2 \, dx,$$

the last term does not depend on $s$, and $\int s(x) U_n'(x) \, dx$ is interpreted as $\int s(x) \, dU_n(x)$. Although not stated explicitly there, the *isotonic inverse estimator* studied in Van Es, Jongbloed and Van Zuijlen (1998) can be interpreted in the same way as the LS estimator considered here. The only difference is that $Q_n$ is minimized over all decreasing rather than convex decreasing functions $[0, \infty)$.

The main reason for considering the survival function $s$ instead of the distribution function $F$ in the definition of the least square estimator is that the survival function is convex and decreasing and, henceforth, we can exploit results from Groeneboom, Jongbloed and Wellner (2001b) more naturally. We now provide conditions on existence of the reciprocal kernel $p$.

LEMMA 2.4. *To each kernel function $k \in \mathcal{K}$ defined in* (3), *there corresponds a reciprocal kernel $p$* (*or "type 1 resolvent"*), *solving the first kind Volterra integral equation of convolution type*

(12) $$(p * k)(x) := \int_0^x p(x - y) k(y) \, dy = x \mathbf{1}_{[0, \infty)}(x).$$

*This function $p$ is increasing, equals zero on $(-\infty, 0)$ and satisfies $p(0+) = 1/k(0+)$. Moreover, $\lim_{t \to \infty} t^{-1} p(t) = 1$. If, in addition, $k$ is smooth in the sense that it can be written as*

(13) $$k(x) = k(0+) - \int_0^x \kappa(y) \, dy = \int_x^\infty \kappa(y) \, dy,$$

*for a Lipschitz continuous nonnegative function $\kappa$ on $(0, \infty)$, then the function $p$ admits a representation*

$$p(t) = \frac{1}{k(0+)} + L(t) = \frac{1}{k(0+)} + \int_0^t \ell(s)\, ds \tag{14}$$

*for a nonnegative continuous function $\ell$ on $(0, \infty)$ that is Lipschitz continuous on each bounded interval.*

REMARK 2.5. For some kernels $k \in \mathcal{K}$, $p$ is explicitly known. For example, $p(t) = (1+t)\mathbf{1}_{[0,\infty)}(t)$ for the standard exponential $k$ and $p(t) = (1 + \lfloor t \rfloor)\mathbf{1}_{[0,\infty)}(t)$ for the uniform$(0,1)$ kernel $k$. For other situations $p$ can be easily approximated numerically using numerical integration procedures.

PROOF OF LEMMA 2.4. For the first part we refer to Van Es, Jongbloed and Van Zuijlen (1998) and Pipkin (1991), Chapter 6. For the result on smooth kernels, consider the Volterra convolution integral equation of the second kind

$$\ell(t) - \int_0^t \frac{\kappa(t-u)}{k(0+)}\ell(u)\, du = \frac{\kappa(t)}{k(0+)^2} \tag{15}$$

and note that if $\ell$ solves this equation, $p$ defined in (14) solves (12). Existence of a continuous solution to (15) is guaranteed by Theorem 3.5 in Gripenberg, Londen and Staffans (1990) because $\kappa$ is continuous. Using Lipschitz continuity of $\kappa$, Lipschitz continuity of $\ell$ follows. Indeed, denote the Lipschitz constant of $\kappa$ by $K$, and let $t \in [0, M]$ and $h > 0$ sufficiently small. Then

$$|\ell(t+h) - \ell(t)| = \left| \int_0^t \frac{\kappa(t-u+h) - \kappa(t-u)}{k(0+)} \ell(u)\, du \right.$$
$$\left. + \int_t^{t+h} \frac{\kappa(t+h-u)}{k(0+)} \ell(u)\, du + \frac{\kappa(t+h) - \kappa(t)}{k(0+)^2} \right|$$
$$\leq \left\{ \frac{K}{k(0+)} \sup_{[0,M]} |\ell(u)| \left(1 + \sup_{[0,M]} |\kappa(u)|\right) + \frac{K}{k(0+)^2} \right\} h = c_M h.$$

The result now follows from continuity of both $\ell$ and $\kappa$ on the compact interval $[0, M]$. □

ASSUMPTION 2.6. Throughout the rest of the paper we will assume that $k$ admits representation (13) with Lipschitz continuous nonnegative function $\kappa$.

REMARK 2.7. Note that $U_n$ is a right-continuous function. The limit behavior of $p$ implies that $U_n(x) = o(x)$, as $x \to \infty$. It is obvious that $U_n(x) = x$ for $x \in [0, Z_{(1)})$ and that $U_n$ has negative jumps of size $\frac{1}{n}p(0)$ at all observation points.

There are two natural ways to define the least squares estimator. The first is to define it as the minimizer of $Q_n$ over the set $\mathcal{S}$, as done above. A drawback of this approach is that additional assumptions on $k$ are needed to show that the estimator $\hat{s}_n$ is well defined and to derive its asymptotic properties. We follow an alternative approach (avoiding these conditions) where we define the least squares estimator as the minimizer of $Q_n$ over the set

(16)
$$\mathcal{S}_n = \{s : s \text{ convex and decreasing,}$$
$$s(0) = 1, s(Z_{(n)}) = 0, s \text{ piecewise linear with kinks only in } \mathcal{Z}_n\}.$$

THEOREM 2.8. *The least squares estimator $\tilde{s}_n$, defined as the minimizer of $Q_n$ over $\mathcal{S}_n$, exists uniquely.*

PROOF. Uniqueness is immediate from strict convexity of $Q_n$. For existence, note that any $s \in \mathcal{S}_n$ can be written as $s = \sum_{i=1}^n \alpha_i s_{Z_i}$, where $s_\theta = 1 - F_\theta$ [with $F_\theta$ defined in (6)], all $\alpha_i \in [0,1]$ and $\sum_{i=1}^n \alpha_i = 1$. Hence, the minimization problem is equivalent to that of minimizing

$$(\alpha_1, \ldots, \alpha_n) \mapsto \frac{1}{2} \sum_{i=1}^n \sum_{j=1}^n \alpha_i \alpha_j \int s_{Z_i} s_{Z_j} \, dx - \sum_{i=1}^n \alpha_i \int s_{Z_i} \, dU_n$$

over the set $\mathcal{C} = \{\alpha_i \in [0,1] \, (i=1, \ldots, n), \sum_{i=1}^n \alpha_i = 1\}$. The existence now follows from the compactness of $\mathcal{C}$ and the continuity of the mapping in the preceding display. □

REMARK 2.9. The following argument shows why we can restrict the minimization to functions that equal one at zero. To show that $\hat{s}_n(0) = 1$, note that the integral in objective function (11) can be split in the regions $[0, Z_{(1)})$ and $[Z_{(1)}, Z_{(n)}]$. The first part is $\frac{1}{2} \int_0^{Z_{(1)}} s(x)(s(x) - 2) \, dx$, where the convex integrand is minimized pointwisely by taking $s(x) = 1$. Hence, for any $s \in \mathcal{S}$ with $s(0) < 1$, the objective function can be decreased by moving $s$ on $[0, Z_{(1)})$ as closely as possible to one. This boils down to changing it to the linear function connecting $(0, 1)$ with $(Z_{(1)}, s(Z_{(1)}))$.

We now state necessary and sufficient conditions that characterize $\tilde{s}_n$.

THEOREM 2.10. *The function $s$ minimizes $Q_n$ over all functions in $\mathcal{S}_n$, if and only if for all $\theta \in \mathcal{Z}_n$*

$$H_n(\theta; s) = \int_{t=0}^\theta \int_{v=0}^t s(v) \, dv \, dt - \theta \left( \int_0^\infty s(t)^2 \, dt - \int_0^\infty s(t) \, dU_n(t) \right)$$

(17)
$$\geq \int_0^\theta U_n(t)\,dt = Y_n(\theta),$$

*with equality whenever $\theta$ is a kink of $s$.*

PROOF. For necessity, assume $s$ minimizes $Q_n$ over $\mathcal{S}_n$. Because $s + \varepsilon(s_\theta - s) \in \mathcal{S}_n$ for all $\theta \in \mathcal{Z}_n$ and $\varepsilon \in [0,1]$, and $s$ minimizes $Q_n$ over $\mathcal{S}_n$, we have that
$$\lim_{\varepsilon \downarrow 0} \varepsilon^{-1}(Q_n(s + \varepsilon(s_\theta - s)) - Q_n(s)) \geq 0.$$

Writing out this limit, we get
$$\int_0^\infty s(x)(s_\theta(x) - s(x))\,dx - \int_0^\infty (s_\theta(x) - s(x))\,dU_n(x) \geq 0 \qquad \forall \theta \in \mathcal{Z}_n.$$

Denote, for the moment, by $\bar{s}$ the primitive of $s$, which is zero at zero. Then we have
$$\int_0^\infty s(x)s_\theta(x)\,dx = \int_0^\theta s_\theta(x)\,d\bar{s}(x) = \frac{1}{\theta}\int_0^\theta \bar{s}(x)\,dx$$

and
$$\int_0^\infty s_\theta(x)\,dU_n(x) = \frac{1}{\theta}\int_0^\theta U_n(x)\,dx.$$

This leads to the necessary inequality for optimality given in (17).

Now, for sufficiency, suppose $s$ satisfies conditions (17). Let $\tilde{s} = \int s_\theta\,d\tilde{\mu}(\theta) \in \mathcal{S}_n$, arbitrary. Define the function $\varepsilon \mapsto \varphi(\varepsilon) := Q_n(s + \varepsilon(\tilde{s} - s))$, which is convex on $[0,1]$. Moreover, $Q_n(\tilde{s}) = \varphi(1) \geq \varphi(0) + \varphi'(0) = Q_n(s) + \varphi'(0)$, where the derivative is interpreted as right derivative. Hence, $s$ minimizes $Q_n$ over $\mathcal{S}_n$ if $\varphi'(0) \geq 0$. To see that this holds, note that
$$\varphi'(0) = \int_{\theta > 0} \frac{1}{\theta}(H_n(\theta;s) - Y_n(\theta))\,d\tilde{\mu}(\theta) \geq 0.$$

If we take $\tilde{s} = s$, then we obtain an equality in this display. This implies that, for all $\theta$ where $s$ has a kink, $H_n(\theta;s) = Y_n(\theta)$. $\square$

Figures 1 and 2 show the maximum likelihood estimator and least squares estimator for the case that the "true" distribution function $F$ equals $F(x) = \min(\sqrt{x/5}, 1)$ $(x > 0)$. In Figure 1 the noise is standard exponentially distributed, and in Figure 2 the noise is sampled from the distribution with density $k(x) = 2(1-x)\mathbf{1}_{[0,1]}(x)$. The sample sizes were taken equal to 10 and 100. The estimators were calculated using the algorithms described in Section 4. Figure 3 gives a plot corresponding to the left-hand side picture of Figure 1. It shows that the MLE and LSE satisfy the characterizations of Theorems 2.3 and 2.10, respectively.

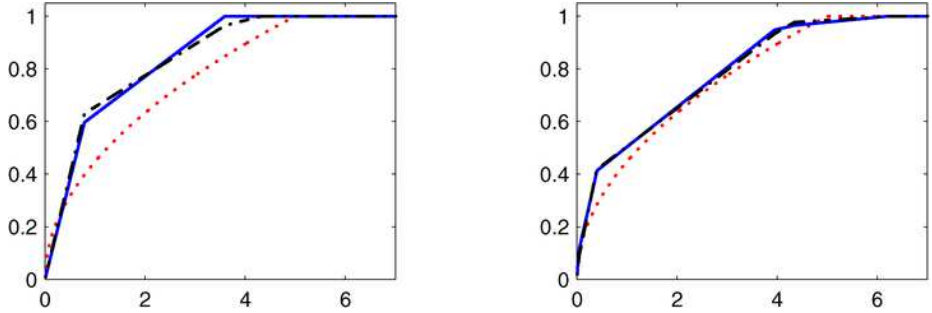

FIG. 1. *Deconvolution with $k(x) = e^{-x}\mathbf{1}_{[0,\infty)}(x)$. Left: $n = 10$. Right: $n = 100$. True: red dotted; MLE: blue solid; LSE: black dash-dotted.*

**3. Consistency of the estimators.** In Theorems 3.1 and 3.3 we prove consistency of the maximum likelihood and least squares estimators, respectively.

3.1. *Maximum likelihood.*

THEOREM 3.1. *Let $k \in \mathcal{K}$ satisfy Assumption 2.6. Then, almost surely, $\|\hat{F}_n - F_0\|_\infty \to 0$. That is, the MLE is strongly uniformly consistent. In addition, we have for all $x > 0$*

(18) $$F_0^l(x) \geq \limsup_{n \to \infty} \hat{F}_n^l(x) \geq \liminf_{n \to \infty} \hat{F}_n^r(x) \geq F_0^r(x).$$

*Here the superscripts "l" and "r" denote left and right derivatives, respectively.*

PROOF. If $\hat{F}_n$ maximizes $l_n$ over $\mathcal{F}$, then, by Theorem 2.3

$$\int \frac{g_{F_0}(z)}{g_{\hat{F}_n}(z)} d\mathbb{G}_n(z) = \int\int \frac{g_\theta(z)}{g_{\hat{F}_n}(z)} d\mu_{F_0}(\theta) d\mathbb{G}_n(z)$$

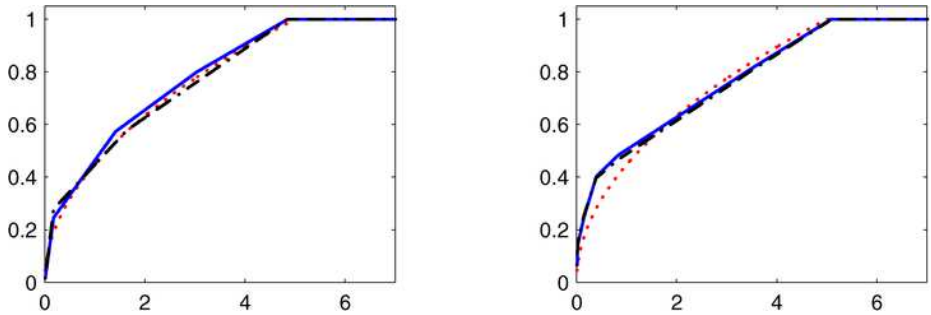

FIG. 2. *Deconvolution with $k(x) = 2(1-x)\mathbf{1}_{[0,1]}(x)$. Left: $n = 10$. Right: $n = 100$. True: red dotted; MLE: blue solid; LSE: black dash-dotted.*

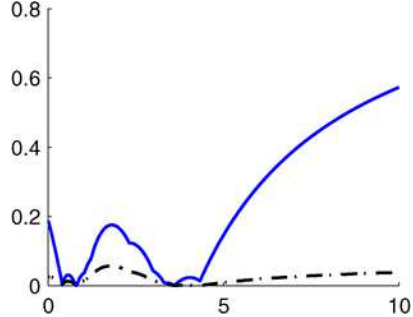

Fig. 3. *Deconvolution with $k(x) = e^{-x}\mathbf{1}_{[0,\infty)}(x)$ ($n = 10$). The curves show that the characterization of the MLE and LSE, as given in Theorems 2.3 and 2.10, respectively, are satisfied. MLE: blue solid; LSE: black dash-dotted.*

(19)
$$= \int \int \frac{g_\theta(z)}{g_{\hat{F}_n}(z)} d\mathbb{G}_n(z) \, d\mu_{F_0}(\theta) \leq 1.$$

By the Glivenko–Cantelli theorem, if $\Omega_0 := \{\|\mathbb{G}_n(\cdot, \omega) - G_0\|_\infty \to 0\}$, where $G_0$ is the distribution function corresponding to $g_{F_0}$, then $\mathbb{P}(\Omega_0) = 1$. Fix $\omega \in \Omega_0$.

Choose an arbitrary subsequence $(m)$ of $(n)$. Using the Helly selection principle, a subsequence $(l)$ of $(m)$ and a concave subdistribution function $\tilde{F}$ on $[0, \infty)$ can be extracted such that $\hat{F}_l(x)$ converges to $\tilde{F}(x)$ for all $x > 0$. By Lemma A.1 in the Appendix, this vague convergence implies for the corresponding convolution densities $\hat{g}_l = g_{\hat{F}_l}$ and (sub) density $\tilde{g} = g_{\tilde{F}}$ that for all closed intervals $I$ in $(0, \infty)$, $\sup_{z \in I} |\hat{g}_l(z) - \tilde{g}(z)| \to 0$ as $l \to \infty$. Following exactly the argument of Theorem 3.2 in Groeneboom, Jongbloed and Wellner (2001b), it can be shown that necessarily $g_0 = \tilde{g}$.

Hence, any subsequence of the sequence $\{\hat{F}_n\}_n$ has a further subsequence $\{\hat{F}_l\}_l$ with $\hat{F}_l \xrightarrow{w} \tilde{F}$ for some $\tilde{F}$. Furthermore, we saw that $\tilde{g} = g_{\tilde{F}} = g_0 = g_{F_0}$. This implies $\tilde{F} = F_0$, so there is only one possible limit for the subsequence. Therefore, for all $\omega \in \Omega_0$, $\hat{F}_n(\omega) \xrightarrow{w} F_0$. Since $F_0$ is concave, it is continuous. This implies that almost surely $\|\hat{F}_n - F_0\|_\infty \to 0$, as $n \to \infty$. The statement in (18) is a consequence of Marshall's lemma [Robertson, Wright and Dykstra (1988), page 332]. □

REMARK 3.2. If we consider the more general setting mentioned in Remark 2.2, then the preceding argument can be extended to prove consistency for this case as well.

3.2. *Least squares.* Next we prove consistency for the least squares estimator. Let $U_0(x) = \int_0^x s_0(y)\,dy$ and define

$$Q_0(s) = \frac{1}{2}\int_0^\infty s(x)^2\,dx - \int_0^\infty s(x)\,dU_0(x)$$
$$= \frac{1}{2}\int_0^\infty (s(x) - s_0(x))^2\,dx - \frac{1}{2}\int_0^\infty s_0(x)^2\,dx.$$

THEOREM 3.3. *Assume $s_0 \in \mathcal{S}$. If we denote the $L^2$-norm of functions on $[0,\infty)$ by $\|\cdot\|_2$, then $\|\tilde{s}_n - s_0\|_2 \xrightarrow{a.s.} 0$ and $\|\tilde{s}_n - s_0\|_\infty \xrightarrow{a.s.} 0$, as $n \to \infty$.*

PROOF. Note that

$$\mathcal{S}_1 \subseteq \mathcal{S}_2 \subseteq \cdots \subseteq \mathcal{S}_n \subseteq \cdots \mathcal{S} \subseteq L^2[0,\infty).$$

For each $i \geq 1$, the set $\mathcal{S}_i$ is closed with respect to the topology induced by the $L^2$-norm. This follows from the fact that $s \in \mathcal{S}_i$ is bounded and piecewise linear, with kinks at at most $i$ points. Furthermore, $\mathcal{S}_i$ is convex. Therefore, the $L^2$-projection of $s_0 \in \mathcal{S}$ onto $\mathcal{S}_i$ exists. Denote the latter by $\Pi_i s_0$. Using the fact that $\tilde{s}_n$ minimizes $Q_n$ over $\mathcal{S}_n$, we get

$$\frac{1}{2}\|\tilde{s}_n - s_0\|_2^2$$
$$= Q_0(\tilde{s}_n) + \frac{1}{2}\|s_0\|_2^2$$
$$= Q_0(\tilde{s}_n) - Q_n(\tilde{s}_n) + Q_n(\tilde{s}_n) + \frac{1}{2}\|s_0\|_2^2$$
$$\leq Q_0(\tilde{s}_n) - Q_n(\tilde{s}_n) + Q_n(\Pi_n s_0) + \frac{1}{2}\|s_0\|_2^2$$
$$= Q_0(\tilde{s}_n) - Q_n(\tilde{s}_n) + Q_n(\Pi_n s_0) - Q_0(\Pi_n s_0) + Q_0(\Pi_n s_0) + \frac{1}{2}\|s_0\|_2^2$$
$$\leq 2\sup_{s \in \mathcal{S}_n}|Q_0(s) - Q_n(s)| + Q_0(\Pi_n s_0) + \frac{1}{2}\|s_0\|_2^2$$
$$= 2\sup_{s \in \mathcal{S}_n}|Q_0(s) - Q_n(s)| + \frac{1}{2}\|\Pi_n s_0 - s_0\|_2^2.$$

On the other hand, we have that

$$U_0(x) - U_n(x) = \int_0^x p(x-y)\,d(\mathbb{G}_n - G_0)(y)$$
$$= \frac{1}{k(0)}(\mathbb{G}_n - G_0)(x) + \int_0^x L(x-y)\,d(\mathbb{G}_n - G_0)(y),$$

where the second equality follows from equation (14). This implies that, for $s \in \mathcal{S}_n$,

$$Q_n(s) - Q_0(s) = \int \left( -\frac{s(x)}{k(0)} + \int_x^\infty s(y)\ell(y-x)\,dy \right) d(\mathbb{G}_n - G_0)(x).$$

Substituting this equality in the preceding inequality gives

$$\|\tilde{s}_n - s_0\|_2^2$$
$$\leq 4 \sup_{s \in \mathcal{S}_n} \left| \int \left( -\frac{s(x)}{k(0)} \right. \right.$$
$$\left. \left. + \int_x^\infty s(y)\ell(y-x)\,dy \right) d(\mathbb{G}_n - G_0)(x) \right| + \|\Pi_n s_0 - s_0\|_2^2$$
$$\leq 4 \sup_{s \in \mathcal{S}} \left| \int \left( -\frac{s(x)}{k(0)} \right. \right.$$
$$\left. \left. + \int_x^\infty s(y)\ell(y-x)\,dy \right) d(\mathbb{G}_n - G_0)(x) \right| + \|\Pi_n s_0 - s_0\|_2^2.$$

Since $\overline{\bigcup_{n=1}^\infty \mathcal{S}_n} = \mathcal{S}$ almost surely, $\|\Pi_n s_0 - s_0\|_2$, tends to zero almost surely, as $n \to \infty$. If the class

$$\left\{ x \mapsto -\frac{s(x)}{k(0)} + \int_x^\infty s(y)\ell(y-x)\,dy,\ s \in \mathcal{S} \right\}$$

is Glivenko–Cantelli, then the first term tends to zero as well. That this class is indeed Glivenko–Cantelli can be seen as follows. First, the class $\mathcal{S}$ itself is Glivenko–Cantelli [Example 3.7.1 in Van de Geer (2000)]. Moreover, $\{v : v(x) = \int_0^\infty s(x+y)\ell(y)\,dy,\ s \in \mathcal{S}\} \subset \mathcal{S}$ is Glivenko–Cantelli for the same reason. Hence, by the triangle inequality, the class consisting of sums of two functions, one from each class, is Glivenko–Cantelli, too.

Now suppose that $\tilde{s}_n$ does not converge to $s_0$ pointwisely. Then there exists a point $x > 0$, and $\varepsilon > 0$ and a subsequence of $n$, such that for all $n$ along this subsequence $|\tilde{s}_n(x) - s_0(x)| > \varepsilon$. Because all $\tilde{s}_n$ and $s_0$ are convex and decreasing, there is a fixed neighborhood of $x$, such that for all $y$ in this neighborhood and $n$ along this subsequence, $|\tilde{s}_n(y) - s_0(y)| > \varepsilon/2$. This implies that $\|\tilde{s}_n - s_0\|_2$ does not converge to zero. Hence, with probability one $\tilde{s}_n(x) \to s_0(x)$ for all $x$, as $n \to \infty$. Uniform consistency follows from this pointwise result because $\tilde{s}_n$ and $s_0$ are convex and decreasing (the proof is similar to the proof of the classical Glivenko–Cantelli theorem). □

**4. Computing the estimators by a support-reduction algorithm.** Both estimators can be computed by the support-reduction algorithm as discussed in Groeneboom, Jongbloed and Wellner (2008). This is an iterative algorithm for minimizing a convex objective function $Q$ over a convex cone or convex hull generated by a parametrized function class. Suppose the objective function is denoted by $Q$, and let the convex cone $\mathcal{F}$ generated by the functions $\{f_\theta : \theta \in \Theta\}$ be given by

$$\mathcal{F} = \left\{ f \Big| f(x) = \int f_\theta(x) \, d\mu(\theta), \mu \text{ is a positive finite measure on } \Theta \right\},$$

where $\Theta$ is some subset of $\mathbb{R}$. (If we minimize over a convex hull, then the measure $\mu$ is a probability measure.) We aim to compute $\hat{f} = \arg\min_{f \in \mathcal{F}} Q(f)$.

Both the computation of the ML estimator and the LS estimator fit within this framework. For the MLE we minimize $Q(f) = -\int \log f(x) \, dG_n(x) + \int f(x) \, dx$ over the convex cone generated by the functions $\{g_\theta : \theta \in \mathcal{Z}_n\}$; for the LSE we minimize $Q(f) = \frac{1}{2} \int f(x)^2 \, dx - \int_0^\infty f(x) \, dU_n(x)$ over the convex hull generated by the functions $\{s_\theta : \theta \in \mathcal{Z}_n\}$. If the solution is given by $\hat{f}_n = \int f_\theta \, d\hat{\mu}_n(\theta)$, then $\hat{F}_n = \int F_\theta \, d\hat{\mu}_n(\theta)$.

The main steps of the algorithm are briefly explained in Section 6.1 of Jongbloed, van der Meulen and van der Vaart (2005). For additional details we refer to Groeneboom, Jongbloed and Wellner (2008). Computational details for the current setup can be found in the Appendix.

**5. Asymptotic lower bound on local minimax risk.** In this section, we derive an asymptotic lower bound to a local minimax risk for estimating the concave distribution function $F_0$ and its (decreasing) derivative $F_0' = f_0$ at an interior point $x_0 > 0$ of its support. On $f$ we impose a local assumption near the point $x_0$:

$$f_0(x) = f_0(x_0) + f_0'(x_0)(x - x_0)(1 + o(1)) \tag{20}$$

as $x \to x_0$ and $f_0'$ is continuous at $x_0$.

Moreover, we assume an integrability condition on $k$ and $F_0$ jointly:

$$\int_{x_0}^\infty \frac{k'(x - x_0)^2}{g_{F_0}(x)} \, dx < \infty. \tag{21}$$

Define for a fixed kernel function $k$ that can be expressed as in (13) the class of sampling densities

$$\mathcal{G} = \Big\{ g : g(z) = \int_0^z k(z - x) f(x) \, dx, \tag{22}$$

$$z \geq 0 \text{ with } f \text{ decreasing density on } (0, \infty) \Big\}.$$

Endow this class of densities with the Hellinger distance,

$$H(g, h) = \left( \frac{1}{2} \int_0^\infty (\sqrt{h(x)} - \sqrt{g(x)})^2 \, dx \right)^{1/2},$$

and denote by $\mathcal{G}_g$ a subset of $\mathcal{G}$ containing a Hellinger ball of positive radius around the fixed density $g \in \mathcal{G}$.

Now consider the problem of estimating the functionals

$$(23) \qquad T_1(g) = F(x_0) \quad \text{and} \quad T_2(g) = f(x_0)$$

based on a sample from density $g$. The difficulty of the problem of estimating a functional $T(g)$ based on a sample of size $n$ from the density $g \in \mathcal{G}$ can be quantified using the concept of a local minimax risk:

$$(24) \qquad R(n, T, \mathcal{G}_g) = \inf_{t_n} \sup_{g \in \mathcal{G}_n} E_{g^{\otimes n}} |t_n(X) - T(g)|,$$

where the infimum is taken over all estimators $t_n$ based on the sample $X = (X_1, \ldots, X_n)$. In Jongbloed (2000), an asymptotic lower bound to this quantity is given in terms of a (local) modulus of continuity $m_g$ of $T$ over $\mathcal{G}_g$:

$$m_g(\varepsilon; T) = \sup\{|T(h) - T(g)| : h \in \mathcal{G}_g \text{ and } H(h, g) \leq \varepsilon\}.$$

In fact, if it can be shown that

$$(25) \qquad m_g(\varepsilon; T) \geq (c\varepsilon)^r (1 + o(1)) \qquad \text{as } \varepsilon \downarrow 0,$$

then [Corollary 2 in Jongbloed (2000)]

$$(26) \qquad \liminf_{n \to \infty} n^{r/2} R(n, T, \mathcal{G}_g) \geq \frac{1}{4} e^{-r/2} \left( \frac{1}{2} c \sqrt{r} \right)^r.$$

THEOREM 5.1. *Let $T_1$ and $T_2$ be defined as in (23) and $\mathcal{G}$ as in (22). Assume that condition (20) is satisfied for the density $f_0$ associated with $g_0$. Then, for the local minimax risk defined in (24), we have*

$$\liminf_{n \to \infty} n^{2/5} R(n, T_1, \mathcal{G}_{g_0}) \geq \frac{1}{8} \left( \frac{|f_0'(x_0)| g_0(x_0)^2}{100 e^2 k(0)^4} \right)^{1/5}$$

*and*

$$\liminf_{n \to \infty} n^{1/5} R(n, T_2, \mathcal{G}_{g_0}) \geq \frac{1}{4} \left( \frac{|f_0'(x_0)|^3 |g_0(x_0)|}{4 e k(0)^2} \right)^{1/5}.$$

PROOF. We construct a family $\{g_\varepsilon : \varepsilon \in [0, \varepsilon_0]\} \subset \mathcal{G}$ with the following properties:

$$(27) \qquad \begin{aligned} |T_1(g_\varepsilon) - T_1(g_0)| &= \tfrac{1}{2} \varepsilon^2 f_0'(x_0)(1 + o(1)) \quad \text{and} \\ |T_2(g_\varepsilon) - T_2(g_0)| &= \varepsilon f_0'(x_0)(1 + o(1)) \end{aligned}$$

for $\varepsilon \downarrow 0$. Moreover,

(28)
$$H(g_\varepsilon, g_0) \leq (c_1 \varepsilon)^{5/2}(1 + o(1))$$
$$\Rightarrow H(g_{\varepsilon^{2/5}/c_1}, g_0) \leq \varepsilon(1 + o(1)) \qquad \text{as } \varepsilon \downarrow 0,$$

where

$$c_1 = \left(\frac{2k(0)^2 f_0'(x_0)^2}{5g_0(x_0)}\right)^{1/5}.$$

This means that for $\varepsilon \downarrow 0$

$$m_g(\varepsilon; T_1) \geq |T_1(g_{\varepsilon^{2/5}/c_1}) - T_1(g_0)|(1 + o(1)) = \frac{|f_0'(x_0)|\varepsilon^{4/5}}{2c_1^2}(1 + o(1))$$

and

$$m_g(\varepsilon; T_2) \geq |T_2(g_{\varepsilon^{2/5}/c_1}) - T_2(g_0)|(1 + o(1)) = \frac{f_0'(x_0)\varepsilon^{2/5}}{c_1}(1 + o(1)).$$

Using these facts in (25) and (26), the statement of the theorem follows.

Let us now define the class $\{g_\varepsilon : \varepsilon \in [0, \varepsilon_0]\}$ and prove (27) and (28). This class is defined based on a perturbation of the underlying distribution function $F_0$. Indeed,

$$g_\varepsilon(z) = \int_0^z k(z - x) \, dF_\varepsilon(x)$$

with

$$F_\varepsilon(x) = \begin{cases} F_0(x), & \text{if } x \notin [x_0 - c_\varepsilon \varepsilon, x_0 + \varepsilon], \\ F_0(x_0 - c_\varepsilon \varepsilon) + (x - x_0 + c_\varepsilon \varepsilon) \\ \qquad \times f_0(x_0 - c_\varepsilon \varepsilon), & \text{if } x \in [x_0 - c_\varepsilon \varepsilon, x_0 - \varepsilon], \\ F_0(x_0 + \varepsilon) + (x - x_0 - \varepsilon)f_0(x_0 + \varepsilon), & \text{if } x \in (x_0 - \varepsilon, x_0 + \varepsilon]. \end{cases}$$

Here, $c_\varepsilon$ is chosen in such a way that $F_\varepsilon$ is continuous at $x_0 - \varepsilon$. Note that $c_\varepsilon \to 3$ as $\varepsilon \downarrow 0$ and $F_\varepsilon$ is a concave distribution function on $[0, \infty)$, for all small values of $\varepsilon$. By assumption (20), the statements in (27) follow immediately. A proof of (28) is given in the Appendix. □

**6. Asymptotic distribution theory for the LS-estimator.** Theorem 2.10 gives a characterization of the least squares estimator that can be used to derive the limit behavior of the estimator at a fixed point. Let $T_n \subset \mathcal{Z}_n = \{Z_1, \ldots, Z_n\}$ denote the set of bend points of $\tilde{s}_n$.

In this section we prove the following result.

THEOREM 6.1. *Suppose that $s_0$ is twice continuously differentiable in a neighborhood of $x_0$, with strictly positive second derivative. Then,*

$$\begin{pmatrix} n^{2/5} c_1(s_0, k)(\tilde{s}_n(x_0) - s_0(x_0)) \\ n^{1/5} c_2(s_0, k)(\tilde{s}'_n(x_0) - s'_0(x_0)) \end{pmatrix} \xrightarrow{d} \begin{pmatrix} H''(0) \\ H'''(0) \end{pmatrix}. \tag{29}$$

*Here $(H''(0), H'''(0))$ are the second and third derivatives at zero of the invelope $H$ of the stochastic process*

$$Y(t) = \int_0^t W(s)\, ds + t^4$$

*(where $W$ is standard two-sided Brownian motion), introduced in Theorem 2.1 of Groeneboom, Jongbloed and Wellner (2001a). The constants $c_1$ and $c_2$ are given by*

$$c_1(s_0, k) = \left( \frac{24 k(0)^4}{g_0(x_0)^2 s''_0(x_0)} \right)^{1/5} \quad and \quad c_2(s_0, k) = \left( \frac{24}{s''_0(x_0)} \right)^{3/5} \left( \frac{k(0)^2}{g_0(x_0)} \right)^{1/5}.$$

PROOF. Consider the processes

$$H_n(x) = \int_0^x \int_0^y \tilde{s}_n(u)\, du\, dy - x \left( \int \tilde{s}_n(u)^2\, du - \int \tilde{s}_n(u)\, dU_n(u) \right) \tag{30}$$

and

$$Y_n(x) = \int_0^x U_n(y)\, dy.$$

By Theorem 2.10, the characterization of the LS estimator can be written as

$$Y_n(x) \begin{cases} \leq H_n(x), & \text{for all } x \in \mathcal{Z}_n, \\ = H_n(x), & \text{for all } x \in T_n. \end{cases}$$

Now define, for $t \in [-n^{1/5} x_0, \infty)$, localized versions of $Y_n$ and $H_n$:

$$Y_n^{\text{loc}}(t) = n^{4/5}(Y_n(x_0 + n^{-1/5} t)$$
$$- Y_n(x_0) - n^{-1/5} t Y'_n(x_0) - \frac{1}{2} n^{-2/5} t^2 s_0(x_0) - \frac{1}{6} n^{-3/5} t^2 s'_0(x_0))$$
$$= n^{4/5} \int_{x_0}^{x_0 + n^{-1/5} t} \left( U_n(v) - U_n(x_0) \right.$$
$$\left. - \int_{x_0}^v (s_0(x_0) + (u - x_0) s'_0(x_0))\, du \right) dv$$

and

$$H_n^{\text{loc}}(t) = n^{4/5}(H_n(x_0 + n^{-1/5} t) - H_n(x_0)$$
$$- n^{-1/5} t H'_n(x_0) - \tfrac{1}{2} n^{-2/5} t^2 s_0(x_0) - \tfrac{1}{6} n^{-3/5} t^3 s'_0(x_0))$$
$$+ A_n + B_n t,$$

where

(31) $\quad A_n = n^{4/5}(H_n(x_0) - Y_n(x_0))\quad \text{and}\quad B_n = n^{3/5}(H'_n(x_0) - Y'_n(x_0)).$

By Lemma A.2, the random variables $A_n$ and $B_n$ are tight.

The necessary and sufficient conditions for optimality can then be rewritten as

$$Y_n^{\mathrm{loc}}(t) \begin{cases} \leq H_n^{\mathrm{loc}}(t), & \text{for all } t \in [-n^{1/5}x_0, \infty) \text{ with } x_0 + n^{-1/5}t \in \mathcal{Z}_n, \\ = H_n^{\mathrm{loc}}(t), & \text{for all } t \text{ with } x_0 + n^{-1/5}t \in T_n. \end{cases}$$

If we define the process $Z_n$ by

$$Z_n(t) := n^{3/5}((U_n - U_0)(x_0 + n^{-1/5}t) - (U_n - U_0)(x_0))$$

then the process $Y_n^{\mathrm{loc}}$ can be rewritten as

$$Y_n^{\mathrm{loc}}(t) = n^{4/5} \int_{x_0}^{x_0 + n^{-1/5}t} (U_n(v) - U_n(x_0) - (U_0(v) - U_0(x_0)))\, dv$$

$$+ n^{4/5} \int_{x_0}^{x_0 + n^{-1/5}t} \int_{x_0}^{v} (s_0(u) - s_0(x_0) - (u - x_0)s'_0(x_0))\, du\, dv$$

$$= \int_0^t Z_n(v)\, dv + \frac{1}{24} s''_0(x_0) t^4 + o(1),$$

where for any $c > 0$ the $o(1)$ term is uniformly in $t \in [-c, c]$ as $n$ tends to infinity. By Lemma A.6 and the continuous mapping theorem, it now follows that

$$Y_n^{\mathrm{loc}}(t) \xrightarrow{d} \frac{\sqrt{g(x_0)}}{k(0)} \int_0^t W(s)\, ds + \frac{1}{24} s''_0(x_0) t^4.$$

Now we proceed by rescaling the axes in the necessary conditions for optimality in such a way that the limiting process behavior of $Y_n^{\mathrm{loc}}$ will no longer depend on the underlying functions $s_0$ and $k$. For any $\alpha, \beta > 0$, the necessary and sufficient conditions can be rewritten as

$$\tilde{H}_n^{\mathrm{loc}}(t) := \alpha H_n^{\mathrm{loc}}(\beta t) \begin{cases} \geq \alpha Y_n^{\mathrm{loc}}(\beta t) =: \tilde{Y}_n^{\mathrm{loc}}(t), & \text{for all } t \in [c, c], \\ = \alpha Y_n^{\mathrm{loc}}(\beta t) =: \tilde{Y}_n^{\mathrm{loc}}(t), & \text{for all } t \in [-c, c] \\ & \text{with } x_0 + n^{-1/5}\beta t \in T_n. \end{cases}$$

In the limit, the right-hand side is given by

$$\alpha \frac{\sqrt{g(x_0)}}{k(0)} \int_0^{\beta t} W(s)\, ds + \frac{\alpha \beta^4}{24} s''_0(x_0) t^4.$$

By Brownian scaling, that is, using that for each $\gamma > 0$, $\sqrt{\gamma} W(\cdot/\gamma)$ is Brownian motion whenever $W$ is, we get that in distribution this process is the same as

$$\alpha \beta^{3/2} \frac{\sqrt{g(x_0)}}{k(0)} \int_0^t W(s)\, ds + \frac{\alpha \beta^4}{24} s''_0(x_0) t^4.$$

In order to get a process that does not depend on properties of $g_0$ or $s_0$, we choose $\alpha$ and $\beta$ such that

$$\alpha\beta^{3/2}\frac{\sqrt{g(x_0)}}{k(0)} = 1 \quad \text{and} \quad \frac{\alpha\beta^4}{24}s_0''(x_0) = 1,$$

yielding

$$\alpha = \left(\frac{s_0''(x_0)}{24}\right)^{3/5}\left(\frac{k(0)^2}{g_0(x_0)}\right)^{4/5} \quad \text{and} \quad \beta = \left(\frac{24\sqrt{g(x_0)}}{s_0''(x_0)k(0)}\right)^{2/5}.$$

Note that

$$(\tilde{H}_n^{\text{loc}})''(0) = \alpha\beta^2 n^{2/5}(\tilde{s}_n(x_0) - s_0(x_0)) = c_1(s_0, k)n^{2/5}(\tilde{s}_n(x_0) - s_0(x_0))$$

and

$$(\tilde{H}_n^{\text{loc}})'''(0) = \alpha\beta^3 n^{1/5}(\tilde{s}_n'(x_0) - s_0'(x_0)) = c_2(s_0, k)n^{1/5}(\tilde{s}_n'(x_0) - s_0'(x_0)).$$

From this point on, essentially the same reasoning can be followed as in the proof of Theorem 6.3 in Groeneboom Jongbloed and Wellner (2001b). Indeed, the necessary and sufficient conditions for optimality can be pushed to the limiting characterization related to the process studied in [Groeneboom, Jongbloed and Wellner (2001b), pages 1689–1690], where also Lemma A.4 is needed to use their tightness argument. This leads to the convergence of the vector $((\tilde{H}_n^{\text{loc}})'''(0), (\tilde{H}_n^{\text{loc}})''(0))$, as described in (29). □

REMARK 6.2. Because $s_0' = -f_0$ by definition, the asymptotic standard deviations of $\tilde{s}_n$ and $\tilde{s}_n'$ coincide with the asymptotic bounds on the minimax risk given in Theorem 5.1, apart from some constants not depending on the underlying functions $s_0$ and $k$.

## APPENDIX

LEMMA A.1. *Let $F_n$ be a sequence of concave distribution functions on $[0,\infty)$ converging to the concave (sub)distribution function $F$ pointwisely on $(0,\infty)$ (i.e., the corresponding sequence of distributions converges vaguely to the subdistribution corresponding to $F$). Let $k$ be a density on $(0,\infty)$ satisfying Assumption 2.6. Denote by $g_n$ and $g$ the convolutions of $k$ with $F_n$ and $F$ respectively. Then, $g_n$ converges to $g$ uniformly on closed bounded intervals not containing 0.*

PROOF. Denote for $p = 1, 2, \ldots$ by $k^{(p)}$ compactly supported functions such that for all $p$, $0 \leq k^{(p)} \leq k$ and such that $\|k - k^{(p)}\|_1 \leq 1/p$. Choose arbitrary $M > 0$, and define $\|g\|_{1,M} = \int_0^M |g(z)|\,dz$ by the triangle inequality

$$(32) \qquad \|g_n - g\|_{1,M} \leq \|g_n - g_n^{(p)}\|_1 + \|g_n^{(p)} - g^{(p)}\|_{1,M} + \|g - g^{(p)}\|_1,$$

where $g_n^{(p)} = k^{(p)} * dF_n$ and $g^{(p)} = k^{(p)} * dF$. Now, choose $\varepsilon > 0$ and take $p > 3/\varepsilon$. For the last term in (32) we have, via Fubini,

$$\|g - g^{(p)}\|_1 = \int_0^\infty \int_0^z (k(z-x) - k^{(p)}(z-x))\, dF(x)\, dz$$
$$\leq \|k - k^{(p)}\|_1 \leq 1/p < \varepsilon/3.$$

The first term in (32) is smaller than $\varepsilon/3$ for the same reason. By the assumed vague convergence, we have for all $z$, $|g_n^{(p)}(z) - g^{(p)}(z)| \to 0$ because $k^{(p)}$ is bounded, continuous and has bounded support. Because $g^{(p)}(z) \leq g(z) \leq k(0+)$, $\|g_n^{(p)} - g^{(p)}\|_{1,M} < \varepsilon/3$ for $n$ sufficiently large by dominated convergence. Now, consider for $\eta > 1$ an interval $[1/\eta, \eta]$. Note that on this interval the densities of $F_n$ and $F$ necessarily take values in the interval $[0, \eta]$. This means that all $g_n$ and $g$ are Lipschitz continuous with constant $\|\kappa\|_\infty + k(0)\eta$:

$$|g(z+h) - g(z)| \leq \int_0^z |k(z+h-x) - k(z-x)|\, dF(x)$$
$$+ \int_z^{z+h} k(z+h-x) f(x)\, dx$$
$$\leq h(\|\kappa\|_\infty + k(0)\eta).$$

This, together with the $\|\cdot\|_{1,M}$ convergence, implies the uniform convergence on $[1/\eta, \eta]$. $\square$

**Computational details for the maximum likelihood estimator.** We aim to minimize

$$l_n(g) = -\int \log g(x)\, d\mathbb{G}_n(x) + \int g(x)\, dx$$

over the set

$$\mathcal{G} := \left\{ g : g(x) = \int_{[0,\infty)} g_\theta(x)\, d\mu(\theta), \mu \text{ is a positive finite measure} \right\}.$$

The addition of the $\int g(x)\, dx$-term in the objective function enables us to minimize over a convex cone instead of a convex hull, since the minimizer of $l_n$ can in fact be shown to be a probability density. By Theorem 2.1, it suffices to consider measures supported on $\mathcal{Z}_n$.

As shown in Section 7 of Groeneboom, Jongbloed and Wellner (2008), given a current iterate $\bar{g}$, instead of $l_n$, we can minimize the local objective function

$$l_n(g; \bar{g}) = \int g(x)\, dx + \int \left\{ \frac{1}{2}\left(\frac{g(x)}{\bar{g}(x)}\right)^2 - 2\frac{g(x)}{\bar{g}(x)} \right\} d\mathbb{G}_n(x),$$

which is a local quadratic approximation of the objective function near $\bar{g}$. This quadratic function can be minimized over the (finitely generated) cone using the support reduction algorithm, yielding

$$\bar{g}_q = \arg\min\{l_n(g;\bar{g}) : g \in \operatorname{cone}(g_\theta : \theta \in \mathcal{Z}_n)\}.$$

The next iterate is then obtained as $\bar{g} + \lambda(\bar{g}_q - \bar{g})$ ($\lambda$ chosen appropriately to assure monotonicity of the algorithm).

We now turn to the details of the support reduction algorithm. To find a new support point (a direction of descent), we first compute

$$l_n(g + \varepsilon g_\theta; \bar{g}) - l_n(g;\bar{g}) = \tfrac{1}{2}\varepsilon^2 c_2(\theta) + \varepsilon c_1(\theta; g).$$

Here,

$$c_1(\theta; g) = 1 - 2\int \frac{g_\theta}{\bar{g}}(x)\,d\mathbb{G}_n(x) + \int \frac{g g_\theta}{\bar{g}^2}(x)\,d\mathbb{G}_n(x),$$

$$c_2(\theta) = \int \frac{g_\theta^2}{\bar{g}^2}\,d\mathbb{G}_n(x).$$

Computations that are completely analogous to those of Section 4 in Groeneboom, Jongbloed and Wellner (2008), then show that the most promising direction is given by

(33) $$\hat\theta = \arg\min_{\theta \in \mathcal{Z}_n} \frac{c_1(\theta;g)}{\sqrt{c_2(\theta)}}.$$

The second step consists of minimizing $l_n(\sum_{i=1}^m \alpha_i g_{\theta_i}; \hat{g})$ over $\alpha_1,\dots,\alpha_m$ (without restrictions on $\alpha_i$). Now

$$l_n\left(\sum_{i=1}^m \alpha_i g_{\theta_i}; \bar{g}\right) = \sum_{i=1}^m \alpha_i \left(1 - 2\int \frac{g_{\theta_i}}{\bar{g}}(x)\,d\mathbb{G}_n(x)\right)$$
$$+ \frac{1}{2}\sum_{i=1}^m \sum_{j=1}^m \alpha_i\alpha_j \int \frac{g_{\theta_i}g_{\theta_j}}{\bar{g}^2}(x)\,d\mathbb{G}_n(x).$$

Differentiating with respect to $\alpha_i$ yields the linear system of equation $A(\alpha_1,\dots,\alpha_m)' = \mathbf{b}$, where

$$A_{i,j} = \int \frac{g_{\theta_i}g_{\theta_j}}{\bar{g}^2}(x)\,d\mathbb{G}_n(x), \qquad \mathbf{b}_i = -1 + 2\int \frac{g_{\theta_i}}{\bar{g}}(x)\,d\mathbb{G}_n(x).$$

**Computational details for the least squares estimator.** The least squares estimator is defined as the minimizer of

$$Q_n(s) = \frac{1}{2}\int_0^\infty s(x)^2\,dx - \int_0^\infty s(x)\,dU_n(x)$$

over the set $\mathcal{S}_n$ as defined in (16). If $s \in \mathcal{S}_n$, then $s(x) = \int_0^\infty s_\theta(x)\,d\mu(\theta)$, where $s_\theta(x) = (1-x/\theta)_+$ and $\mu$ is a probability measure supported on $\mathcal{Z}_n$. In the following, we denote $\langle f, g \rangle = \int f(x)g(x)\,dx$ and $\langle f, dU_n \rangle = \int f(x)\,dU_n(x)$.

In the first step of the support reduction algorithm we look for a direction of descent. Given an iterate $s$, the directional derivative in the direction of $s_\theta$ is given by

$$c_1(\theta; s) = \lim_{\varepsilon \downarrow 0} \varepsilon^{-1}(Q_n(s + \varepsilon s_\theta) - Q_n(s)) = \langle s, s_\theta \rangle - \langle s_\theta, dU_n \rangle.$$

The new support point is given by $\hat{\theta} = \arg\min_{\theta \in \mathcal{Z}_n} c_1(\theta; s)$. By Theorem 2.10, the optimal solution $\hat{s}$ satisfies $c_1(\theta; \hat{s}) \geq \langle \hat{s}, \hat{s} \rangle - \langle \hat{s}, dU_n \rangle$.

The second step of the algorithm consists of minimizing $Q_n(\sum_{i=1}^m \alpha_i s_{\theta_i})$ over all $\alpha_i$, such that $\sum_{i=1}^m \alpha_i = 1$. If $m = 1$, we simply have $\alpha_1 = 1$. Else, we set $\alpha_1 = 1 - \sum_{i=2}^m \alpha_i$ and minimize over $\alpha_2, \ldots, \alpha_m$ (without restrictions). We can write

$$Q_n\left(\sum_{i=1}^m \alpha_i s_{\theta_i}\right) = \frac{1}{2} \sum_{i=1}^m \sum_{j=1}^m \alpha_i \alpha_j \langle s_{\theta_i}, s_{\theta_j} \rangle - \sum_{i=1}^m \alpha_i \langle s_{\theta_i}, dU_n \rangle$$

$$= \frac{1}{2} \alpha_1^2 \langle s_{\theta_1}, s_{\theta_1} \rangle + \alpha_1 \sum_{i=2}^m \alpha_i \langle s_{\theta_1}, s_{\theta_i} \rangle + \frac{1}{2} \sum_{i=2}^m \sum_{j=2}^m \alpha_i \alpha_j \langle s_{\theta_i}, s_{\theta_j} \rangle$$

$$- \alpha_1 \langle s_{\theta_1}, dU_n \rangle - \sum_{i=2}^m \alpha_i \langle s_{\theta_i}, dU_n \rangle.$$

Differentiating with respect to $\alpha_i$ $(i = 2, \ldots, m)$, yields the linear system of equations $A(\alpha_2, \ldots, \alpha_m)' = \mathbf{b}$, where

$$A_{i-1, j-1} = \langle s_{\theta_1} - s_{\theta_i}, s_{\theta_1} - s_{\theta_j} \rangle, \qquad i, j = 2, \ldots, m,$$

and

$$\mathbf{b}_{i-1} = \langle s_{\theta_1} - s_{\theta_i}, s_{\theta_1} - dU_n \rangle, \qquad i = 2, \ldots, m.$$

**Proof of (28).** For ease of notation we shall omit subscripts on $f$ and $g$ in the proof. Thus, we write $f$ instead of $f_0$. We use Lemma 2 from Jongbloed (2000), which states that

$$H^2(g_\varepsilon, g) \sim \frac{1}{8} \int_{\{x\,:\,g(x) > 0\}} \frac{(g_\varepsilon(x) - g(x))^2}{g(x)}\,dx = I_\varepsilon^{(1)} + I_\varepsilon^{(2)} + I_\varepsilon^{(3)} \qquad \text{as } \varepsilon \downarrow 0,$$

where $I_\varepsilon^{(1)}$, $I_\varepsilon^{(2)}$ and $I_\varepsilon^{(3)}$ are defined as the integral over the regions $[x_0 - c_\varepsilon \varepsilon, x_0 - \varepsilon]$, $(x_0 - \varepsilon, x_0 + \varepsilon]$ and $(x_0 + \varepsilon, \infty)$ respectively. Note that, for all $x \geq 0$,

$$g(x) - g_\varepsilon(x) = \int_{x_0 - c_\varepsilon \varepsilon}^{x_0 - \varepsilon} k(x - u)(f(u) - f(x_0 - c_\varepsilon \varepsilon))\,du$$

$$+ \int_{x_0 - \varepsilon}^{x_0 + \varepsilon} k(x - u)(f(u) - f(x_0 + \varepsilon))\,du$$

and that, for $x < x_0 - c_\varepsilon \varepsilon$, this difference is zero, since $k(x) = 0$ for $x < 0$. For $x \in [x_0 - c_\varepsilon \varepsilon, x_0 - \varepsilon]$, we have that $g(x) - g_\varepsilon(x) = \int_{x_0 - c_\varepsilon \varepsilon}^{x} k(x - u)(f(u) - f(x_0 - c_\varepsilon \varepsilon)) \, du$. Since $k$ satisfies (13), $\sup_{u \in (x_0 - c_\varepsilon \varepsilon, x)} |k(x - u) - k(0)| = o(1)$ as $\varepsilon \downarrow 0$. Furthermore, condition (20) implies

$$f(u) - f(x_0 - c_\varepsilon \varepsilon) = (u - x_0 + c_\varepsilon \varepsilon) f'(\xi),$$
$$\xi \in (x_0 - c_\varepsilon \varepsilon, u) \subseteq (x_0 - c_\varepsilon \varepsilon, x_0 - \varepsilon).$$

If $\varepsilon \downarrow 0$, then $\xi \to x_0$ and $f'(\xi) \to f'(x_0)$, since $f'$ is continuous at $x_0$. Hence,

$$g(x) - g_\varepsilon(x) = \int_{x_0 - c_\varepsilon \varepsilon}^{x} (k(0) + o(1))(u - x_0 + c_\varepsilon \varepsilon)(f'(x_0) + o(1)) \, du$$

(34)
$$= \frac{1}{2} k(0) f'(x_0) [(u - x_0 + c_\varepsilon \varepsilon)^2]_{x_0 - c_\varepsilon \varepsilon}^{x} (1 + o(1))$$

$$= \frac{1}{2} k(0) f'(x_0) (x - x_0 + c_\varepsilon \varepsilon)^2 (1 + o(1)).$$

Hence,

$$I_\varepsilon^{(1)} = \frac{1}{8} \int_{x_0 - c_\varepsilon \varepsilon}^{x_0 - \varepsilon} \frac{(g_\varepsilon(x) - g(x))^2}{g(x)} \, dx$$

$$= \frac{k(0)^2 f'(x_0)^2}{32} \int_{x_0 - c_\varepsilon \varepsilon}^{x_0 - \varepsilon} \frac{(x - x_0 + c_\varepsilon \varepsilon)^4}{g(x)} (1 + o(1)) \, dx$$

$$= \frac{k(0)^2 f'(x_0)^2}{5 g(x_0)} \varepsilon^5 (1 + o(1)).$$

For $x \in (x_0 - \varepsilon, x_0 + \varepsilon)$,

$$g(x) - g_\varepsilon(x) = \int_{x_0 - c_\varepsilon \varepsilon}^{x_0 - \varepsilon} k(x - u)(f(u) - f(x_0 - c_\varepsilon \varepsilon)) \, du$$
$$+ \int_{x_0 - \varepsilon}^{x} k(x - u)(f(u) - f(x_0 + \varepsilon)) \, du.$$

In exactly the same manner as the previous case, we can find asymptotic order relations for this expression. For the first term we get, from (34),

$$\int_{x_0 - c_\varepsilon \varepsilon}^{x_0 - \varepsilon} k(x - u)(f(u) - f(x_0 - c_\varepsilon \varepsilon)) \, du = 2k(0) f'(x_0) \varepsilon^2 (1 + o(1)).$$

For the second term we get

$$\int_{x_0 - \varepsilon}^{x} k(x - u)(f(u) - f(x_0 + \varepsilon)) \, du$$
$$= \frac{1}{2} k(0) f'(x_0) [(x - x_0 - \varepsilon)^2 - 4\varepsilon^2](1 + o(1)).$$

This gives $g(x) - g_\varepsilon(x) = \frac{1}{2}k(0)f'(x_0)(x - x_0 - \varepsilon)^2(1 + o(1))$, and thus

$$I_\varepsilon^{(2)} = \frac{1}{8}\int_{x_0-\varepsilon}^{x_0+\varepsilon} \frac{(g_\varepsilon(x) - g(x))^2}{g(x)}\,dx = \frac{k(0)^2 f'(x_0)^2}{5g(x_0)}\varepsilon^5(1 + o(1)).$$

Now take $x > x_0 + \varepsilon$. Then we can write

$$g(x) - g_\varepsilon(x)$$
$$= \int_{x_0-c_\varepsilon\varepsilon}^{x_0-\varepsilon} k(x-u)(f(u) - f(x_0 - c_\varepsilon\varepsilon))\,du$$
$$+ \int_{x_0-\varepsilon}^{x_0+\varepsilon} k(x-u)(f(u) - f(x_0 + \varepsilon))\,du$$
$$= \int_{x_0-\varepsilon}^{x_0+\varepsilon} \{k(x-u)[f(u) - f(x_0 + \varepsilon)]$$
$$+ k(x - u + (c_\varepsilon - 1)\varepsilon)[f(u - (c_\varepsilon - 1)\varepsilon) - f(x_0 - c_\varepsilon\varepsilon)]\}\,du.$$

Next, we use relations like

$$f(u) - f(x_0 + \varepsilon) = (u - x_0 - \varepsilon)f'(x_0)(1 + o(1))$$

and

$$k(x - u) = k(x - x_0) + (x_0 - u)k'(x - x_0)(1 + o(1))$$

to obtain

$$g(x) - g_\varepsilon(x)$$
$$= k'(x - x_0)f'(x_0)\int_{x_0-\varepsilon}^{x_0+\varepsilon} \{(x_0 - u)(u - x_0 - \varepsilon) + \cdots$$
$$+ (u - x_0 + \varepsilon)$$
$$\times (x_0 - u + (c_\varepsilon - 1)\varepsilon)\}\,du\,(1 + o(1))$$
$$= \frac{8}{3}k'(x - x_0)f'(x_0)\varepsilon^3(1 + o(1)).$$

Now

$$I_\varepsilon^{(3)} = \frac{8}{9}f'(x_0)^2\varepsilon^6 \int_{x_0+\varepsilon}^{\infty} \frac{k'(x - x_0)^2}{g(x)}\,dx\,(1 + o(1))$$

by (21)

$$H^2(g_\varepsilon, g) \sim \frac{2k(0)^2 f'(x_0)^2}{5g(x_0)}\varepsilon^5(1 + o(1)) \qquad \text{as } \varepsilon \downarrow 0.$$

**Technical results for deriving the asymptotic distribution.** In what follows we assume, as in Theorem 6.1, that $s_0$ is twice continuously differentiable in a neighborhood of $x_0$, with strictly positive second derivative.

LEMMA A.2. *The random variables $A_n$ and $B_n$ as defined in* (31) *are tight.*

To be able to prove the lemma, we first need to prove several other lemmas.

**Distance between successive bend points of the estimator.** Recall that $T_n$ denotes the set of bend-points of $\tilde{s}_n$. For a sequence $\xi_n$ converging to $x_0$, define the bend points to the left and right of $\xi_n$ by

$$(35) \quad \tau_n^- = \max\{x \in T_n : x \leq \xi_n\} \quad \text{and} \quad \tau_n^+ = \min\{x \in T_n : x > \xi_n\}.$$

By consistency and the local assumption of strict convexity of $s_0$ in a neighborhood of $x_0$, it follows that $\tau_n^+ - \tau_n^- \xrightarrow{P} 0$ as $n \to \infty$. The lemma below strengthens this to a rate result for $\tau_n^+ - \tau_n^-$ that is used to obtain a rate result for the LS estimator itself.

LEMMA A.3. *Let $\xi_n$ be a sequence converging to $x_0$. Let $\tau_n^+$ and $\tau_n^-$ be defined according to* (35) *Then,*

$$\tau_n^+ - \tau_n^- = O_P(n^{-1/5}).$$

PROOF. Define, for $u < v$, the v-shaped functions connecting the points $(u, 1), ((u+v)/2, -1)$, and $(v, 1)$, also used in Mammen (1991):

$$f_{u,v}(x) = \left(\frac{4}{v-u}\left|x - \frac{u+v}{2}\right| - 1\right)1_{[u,v]}(x).$$

Note that

$$(36) \quad \int f_{u,v}(x)\,dx = \int x f_{u,v}(x)\,dx = 0 \quad \text{and}$$

$$\int x^2 f_{u,v}(x)\,dx = (v-u)^3/24.$$

Now, take $u = \tau_n^-$ and $v = \tau_n^+$ and define the function $\tilde{f}_{u,v}$ as follows. First, set $\tilde{f}_{u,v}(0) = 0$. For $x = Z_1, \ldots, Z_n$, let $\tilde{f}_{u,v}(x) := f_{u,v}(Z_i)$. In between these points define $\tilde{f}_{u,v}$ by linear interpolation. For $x > Z_{(n)}$, $\tilde{f}_{u,v}(x) = 0$. Note that $\tilde{f}_{u,v}$ and $f_{u,v}$ only differ on the spacings containing $u$, $(u+v)/2$ and $v$.

Using (36) and that the maximal distance between successive order statistics is $O_P(n^{-1} \log n)$, it follows that

$$\text{(37)} \qquad \int \tilde{f}_{u,v}(x)\,dx = O_P\left(\frac{\log n}{n}\right), \qquad \int x\tilde{f}_{u,v}(x)\,dx = O_P\left(\frac{\log n}{n}\right)$$

and

$$\text{(38)} \qquad \int x^2 \tilde{f}_{u,v}(x)\,dx = (v-u)^3/24 + O_P\left(\frac{\log n}{n}\right).$$

Observe that, for small positive $\varepsilon$, the function $\hat{s}_n + \varepsilon \tilde{f}_{u,v} \in \mathcal{S}_n$. This implies that

$$\lim_{\varepsilon \downarrow 0} \varepsilon^{-1}(Q(\hat{s}_n + \varepsilon \tilde{f}_{u,v}) - Q(\hat{s}_n)) \geq 0$$

hence $\int \hat{s}_n(x)\tilde{f}_{u,v}(x)\,dx - \int \tilde{f}_{u,v}(x)\,dU_n(x) \geq 0$.

Note that, by (37) and the fact that $\hat{s}_n$ is linear on $[u,v]$, the first term is $O_P(n^{-1}\log n)$. Hence,

$$\text{(39)} \qquad \int \tilde{f}_{u,v}(x)\,d(U_n - U_0)(x) + \int \tilde{f}_{u,v}(x)\,dU_0(x) \leq O_P\left(\frac{\log n}{n}\right).$$

Using that $U_0' = s_0$ and using a Taylor expansion for $s_0$ as well as (37) and (38), we can write for the second term in (39)

$$\int \tilde{f}_{u,v}(x)\,dU_0(x) = \frac{1}{48} s_0''(x_0)(v-u)^3 + O_P\left(\frac{\log n}{n}\right) + o((v-u)^3)$$

yielding

$$\text{(40)} \qquad \begin{aligned}\int \tilde{f}_{u,v}(x)\,d(U_n - U_0)(x) &+ \frac{1}{48} s_0''(x_0)(v-u)^3 \\ &\leq O_P\left(\frac{\log n}{n}\right) + o((v-u)^3).\end{aligned}$$

For the first term in (40), we have

$$\int \tilde{f}_{u,v}(x)\,d(U_n - U_0)(x) = (U_n - U_0)(v) - (U_n - U_0)(u)$$

$$+ \frac{4}{v-u}\left\{\int_u^{(u+v)/2} - \int_{(u+v)/2}^v\right\}(U_n - U_0)(x)\,dx$$

$$= \int \varphi_{u,v}(x)\,d(\mathbb{G}_n - G_0)(x) + O_P\left(\frac{\log n}{n}\right),$$

using the notation $\bar{p}(x) = \int_0^x p(y)\,dy$,

$$\varphi_{u,v}(x) = p(u-x) - p(v-x)$$
$$- \frac{4}{v-u}\left(\bar{p}(u-x) - 2\bar{p}\left(\frac{u+v}{2} - x\right) + \bar{p}(v-x)\right).$$

We now show that, for any $\varepsilon > 0$, by taking $A > 0$ sufficiently large,

$$
\begin{aligned}
(41) \quad P\bigg(&\exists u \in (\xi_n - \delta, \xi_n], v \in (\xi_n, \xi_n + \delta]: \\
&\left|\int \varphi_{u,v}(x) \, d(\mathbb{G}_n - G_0)(x)\right| > \varepsilon(v-u)^3 + An^{-3/5}\bigg)
\end{aligned}
$$

can be made arbitrarily small, uniformly in $n$. To this end, define for $i, j \in K_n = \{1, 2, \ldots, \lceil n^{1/5}\delta \rceil\}$ the sets

$$I_i = (\xi_n - in^{-1/5}, \xi_n - (i-1)n^{-1/5}] \quad \text{and}$$

$$J_j = (\xi_n + (j-1)n^{-1/5}, \xi_n + jn^{-1/5}]$$

and note that the class of functions $\mathcal{F}_{i,j} = \{\varphi_{u,v} : u \in I_i, j \in J_j\}$ is a VC class with envelope

$$
(42) \quad x \mapsto \begin{cases} c(j+i)n^{-1/5}, & \text{for } x \in [0, \xi_n - in^{-1/5}), \\ c, & \text{for } x \in [\xi_n - in^{-1/5}, \xi_n + jn^{-1/5}], \\ 0, & \text{for } x > \xi_n + jn^{-1/5}, \end{cases}
$$

where $c > 0$ is a constant. For deriving this envelope function, we use relation (14) and the Lipschitz continuity of $\ell$. For $y \leq u$,

$$|\varphi_{u,v}(y)| \leq \|\ell\|_\infty (v-u) + \frac{4}{v-u}|p(\xi_{u,v,y})(v-u)/2 - p(\nu_{u,v,y})(v-u)/2|$$

$$\leq \|\ell\|_\infty(v-u) + 2\|\ell\|_\infty|\nu_{u,v,y} - \xi_{u,v,y}| \leq 3\|\ell\|_\infty(v-u).$$

Taking into account that, for $u \in I_i$ and $j \in J_j$, $0 \leq v - u \leq (i+j)n^{-1/5}$, we get the first inequality in (42). The other bounds in (42) can be deduced similarly.

For the probability in (41) we can now write

$$
\begin{aligned}
P\bigg(&\exists i, j \in K_n : \exists u \in I_i, \\
&v \in J_j : \left|\int \varphi_{u,v}(x) \, d(\mathbb{G}_n - G_0)(x)\right| > \varepsilon(v-u)^3 + An^{-3/5}\bigg) \\
\leq P\bigg(&\exists i, j \in K_n : \exists u \in I_i, \\
&v \in J_j : \left|n^{3/5}\int \varphi_{u,v}(x) \, d(\mathbb{G}_n - G_0)(x)\right| > \varepsilon(j+i-2)^3 + A\bigg) \\
\leq P\bigg(&\exists i, j \in K_n : \sup_{u \in I_i, v \in J_j} \left|n^{3/5}\int \varphi_{u,v}(x) \, d(\mathbb{G}_n - G_0)(x)\right| \\
&\qquad\qquad\qquad\qquad\qquad\qquad > \varepsilon(j+i-2)^3 + A\bigg)
\end{aligned}
$$

$$\leq \sum_{i\in K_n} \sum_{j\in K_n} P\bigg(\sup_{u\in I_i, v\in J_j} \bigg|n^{1/10} n^{1/2} \int \varphi_{u,v}(x)\, d(\mathbb{G}_n - G_0)(x)\bigg|$$
$$> \varepsilon(j+i-2)^3 + A\bigg)$$

$$\leq \sum_{i\in K_n} \sum_{j\in K_n} \frac{n^{1/5}}{(\varepsilon(j+i-2)^3 + A)^2}$$
$$\times E\bigg(\sup_{u\in I_i, v\in J_j} \bigg|n^{1/2} \int \varphi_{u,v}(x)\, d(\mathbb{G}_n - G_0)(x)\bigg|\bigg)^2.$$

To bound the expectation in the summand in this expression, we can then use Theorem 2.14.1 in van der Vaart and Wellner (1996), with envelope function (42), yielding, for some positive $c$,

$$E\bigg(\sup_{u\in I_i, v\in J_j} \bigg|n^{1/2} \int \varphi_{u,v}(x)\, d(\mathbb{G}_n - G_0)(x)\bigg|\bigg)^2$$
$$\leq c((i+j)n^{-1/5} + (i+j)^2 n^{-2/5}).$$

This gives, as upper bound for probability (41),

$$\sum_{i=1}^{\infty}\sum_{j=1}^{\infty} \frac{c((i+j)+(i+j)^2 n^{-1/5})}{(\varepsilon(j+i-2)^3 + A)^2} = c\sum_{k=2}^{\infty} \frac{k(k-1)+k^2(k-1)n^{-1/5}}{(\varepsilon(k-2)^3 + A)^2},$$

which, by dominated convergence, can be made arbitrarily small by taking $A$ sufficiently large.

Combining this result with inequality (40), taking $\varepsilon = s_0''(x_0)/96$, we obtain that by taking $A$ sufficiently large, we have with arbitrarily high probability that

$$\frac{n^{3/5} s_0''(x_0)}{48}(\tau_n^+ - \tau_n^-)^3 \leq n^{3/5}\bigg|\int \varphi_{\tau_n^-, \tau_n^+}(x)\, d(\mathbb{G}_n - G_0)(x)\bigg| + O_P\bigg(\frac{\log n}{n^{2/5}}\bigg)$$
$$\leq \frac{n^{3/5} s_0''(x_0)}{96}(\tau_n^+ - \tau_n^-)^3 + A + O_P\bigg(\frac{\log n}{n^{2/5}}\bigg)$$

implying that $\tau_n^+ - \tau_n^- = O_P(n^{-1/5})$. □

**Rate results for the estimator.** The next lemma shows that, in $O_P(n^{-1/5})$ neighborhoods of $x_0$, the minimal value of the difference between $\tilde{s}_n$ and $s_0$ over this neigborhood is $O_P(n^{-2/5})$.

LEMMA A.4. *Let $\xi_n$ be a sequence converging to $x_0$. For any $\varepsilon > 0$ there exist an $M > 1$ and a $c > 0$, such that the following holds with probability* greater than $1 - \varepsilon$. There are bend points $\tau_n^- < \xi_n < \tau_n^+$ of $\tilde{s}_n$ with $n^{1/5} \leq \tau_n^+ - \tau_n^- \leq Mn^{1/5}$, and for any such points we have

$$\inf_{t \in [\tau_n^-, \tau_n^+]} |\tilde{s}_n(t) - s_0(t)| \leq cn^{-2/5} \quad \text{for all } n.$$

PROOF. Applying Lemma A.3 to the sequences $\xi_n \pm n^{-1/5}$ implies that for any $\varepsilon > 0$ we can find an $M > 1$, such that with probability greater than $1 - \varepsilon$ there are bend points of $\tilde{s}_n$ satisfying $\xi_n - Mn^{-1/5} \leq \tau_n^- \leq \xi_n - n^{-1/5} \leq \xi_n + n^{-1/5} \leq \tau_n^+ \leq \xi_n + Mn^{-1/5}$.

Now, fix $\varepsilon > 0$ and define the $M$ and $\tau_n^{\pm}$ accordingly. Define the functions $\varphi_n^{(1)}$ and $\varphi_n^{(2)}$ by

$$\varphi_n^{(1)}(x) = (\tau_n^+ - x)1_{(\tau_n^-, \tau_n^+]}(x) \quad \text{and} \quad \varphi_n^{(2)}(x) = (\tau_n^+ - x)1_{[\tau_n^-, \tau_n^+]}(x)$$

and note that, for $\varepsilon > 0$ sufficiently small, the piecewise linear functions defined by $\tilde{s}_n(z_i) + \varepsilon \varphi_n^{(1)}(z_i)$ and $\tilde{s}_n(z_i) - \varepsilon \varphi_n^{(2)}(z_i)$ (and linear interpolation between observation points) belong to the class $\mathcal{S}_n$. Hence,

$$\lim_{\varepsilon \downarrow 0} \varepsilon^{-1}(Q(\tilde{s}_n + \varepsilon \varphi_n^{(1)}) - Q(\tilde{s}_n)) \geq 0.$$

This implies, taking into account issues related to piecewise linearity of the function via the $O_P(n^{-1} \log n)$ term,

$$(43) \quad \int_{\tau_n^-}^{\tau_n^+} (\tau_n^+ - x)\tilde{s}_n(x)\,dx - \int_{[\tau_n^-, \tau_n^+]} (\tau_n^+ - x)\,dU_n(x) \geq O_P\left(\frac{\log n}{n}\right).$$

Similarly, taking $-\varepsilon \varphi_n^{(2)}$ instead of $\varepsilon \varphi_n^{(1)}$, we obtain

$$(44) \quad \int_{\tau_n^-}^{\tau_n^+} (\tau_n^+ - x)\tilde{s}_n(x)\,dx - \int_{[\tau_n^-, \tau_n^+]} (\tau_n^+ - x)\,dU_n(x) \leq 0.$$

From (43) and (44) we obtain

$$(45) \quad \left| \int_{\tau_n^-}^{\tau_n^+} (\tau_n^+ - x)(\tilde{s}_n(x) - s_0(x))\,dx - \int_{[\tau_n^-, \tau_n^+]} (\tau_n^+ - x)\,d(U_n - U_0)(x) \right|$$
$$= O_P\left(\frac{\log n}{n}\right).$$

Now, suppose that

$$(46) \quad \inf_{x \in [\tau_n^-, \tau_n^+]} |\tilde{s}_n(x) - s_0(x)| > cn^{-2/5}.$$

Then

$$\left| \int_{\tau_n^-}^{\tau_n^+} (\tau_n^+ - x)(\tilde{s}_n(x) - s_0(x))\,dx \right| > c(\tau_n^+ - \tau_n^-)^2 n^{-2/5}$$

which, in view of (45), implies (using that $\tau_n^+ - \tau_n^- \geq n^{-1/5}$)

(47) $\left| \int_{[\tau_n^-, \tau_n^+]} (\tau_n^+ - x) \, d(U_n - U_0)(x) \right| > c(\tau_n^+ - \tau_n^-)^2 n^{-2/5} \geq c n^{-4/5}.$

Also, note that

$$\int_{[\tau_n^-, \tau_n^+]} (\tau_n^+ - x) \, d(U_n - U_0)(x) = \int_{\tau_n^-}^{\tau_n^+} (U_n - U_0)(x) - (U_n - U_0)(\tau_n^-) \, dx$$
$$= O_P(n^{-4/5})$$

by Lemmas A.3 and A.6. Hence, the probability of (46) is smaller than or equal to that of (47), which can be made arbitrarily small by taking $c$ sufficiently large. □

LEMMA A.5. *For each $M > 0$,*

$$\sup_{t \in [-M, M]} |\tilde{s}_n(x_0 + n^{-1/5} t) - s_0(x_0) - n^{-1/5} t s_0'(x_0)| = O_P(n^{-2/5})$$

*and*

$$\sup_{t \in [-M, M]} |\tilde{s}_n'(x_0 + n^{-1/5} t) - s_0'(x_0)| = O_P(n^{-1/5}).$$

PROOF. This follows from Lemmas A.4 and A.3 in the same way Lemma 4.4 follows from Lemmas 4.3 and 4.2 in Groeneboom, Jongbloed and Wellner (2001b). □

PROOF OF LEMMA A.2. Note that the characterization of $\hat{s}_n$ in Theorem 2.10 implies that, for all bend points $\tau_n$ of $\hat{s}_n$,

(48) $H_n(\tau_n) = Y_n(\tau_n) \quad \text{and} \quad H_n'(\tau_n) = Y_n'(\tau_n) + O_P\left(\frac{\log n}{n}\right),$

where the derivative of $Y_n$ is to be interpreted as a right derivative. Choose $\tau_n$, the last bend point of $\hat{s}_n$ before $x_0$. First, consider $B_n$ and observe that

$$B_n = n^{3/5}(H_n'(x_0) - H_n'(\tau_n) + Y_n'(\tau_n) - Y_n'(x_0) + H_n'(\tau_n) - Y_n'(\tau_n))$$
$$= n^{3/5}\left( \int_{\tau_n}^{x_0} \tilde{s}_n(u) \, du - (U_0(x_0) - U_0(\tau_n)) \right)$$
$$- n^{3/5}((U_n - U_0)(x_0) - (U_n - U_0)(\tau_n)) + n^{3/5}(H_n'(\tau_n) - Y_n'(\tau_n)).$$

By (48), the last term is $O_P(n^{-2/5} \log n)$. By Lemmas A.6 and A.3, the second term is $O_P(1)$. To see that the first term is $O_P(1)$ as well, we use a

Taylor expansion of $U_0(x) = \int_0^x s_0(y)\,dy$ around $x_0$,

$$U_0(x_0) - U_0(\tau_n) = (x_0 - \tau_n)U_0'(x_0) - \frac{1}{2}(x_0 - \tau_n)^2 U_0''(x_0)$$

$$+ \frac{1}{6}(x_0 - \tau_n)^3 U_0'''(\xi_n)$$

$$= \int_{\tau_n}^{x_0} (s_0(x_0) + (u - x_0)s_0'(x_0))\,du + \frac{1}{6}(x_0 - \tau_n)^3 s_0''(\xi_n)$$

for $\xi_n \in (\tau_n, x_0)$. Inserting this into the first term gives, for $n$ sufficiently large,

$$n^{3/5} \left| \int_{\tau_n}^{x_0} \tilde{s}_n(u)\,du - (U_0(x_0) - U_0(\tau_n)) \right|$$

$$= n^{3/5} \left| \int_{\tau_n}^{x_0} [\tilde{s}_n(u) - s_0(x_0) - (u - x_0)s_0'(x_0)]\,du - \frac{1}{6} n^{3/5}(x_0 - \tau_n)^3 s_0''(\xi_n) \right|$$

$$\leq n^{3/5}(x_0 - \tau_n) \sup_{u \in [\tau_n, x_0]} |\tilde{s}_n(u) - s_0(x_0) - (u - x_0)s_0'(x_0)|$$

$$+ \frac{1}{3} n^{3/5} s_0''(x_0)(x_0 - \tau_n)^3$$

$$= O_P(1)$$

by Lemmas A.3 and A.5.

Now, for $A_n$ we get

$$A_n = n^{4/5}\{H_n(x_0) - H_n(\tau_n) - (x_0 - \tau_n)H_n'(\tau_n)$$
$$- (Y_n(x_0) - Y_n(\tau_n) - (x_0 - \tau_n)Y_n'(\tau_n))\}$$
$$- n^{4/5}\{(x_0 - \tau_n)(Y_n'(\tau_n) - H_n'(\tau_n))\}.$$

By (48) the second term is $O_P(n^{-1/5}\log n)$. Note that

$$H_n(x_0) - H_n(\tau_n) - (x_0 - \tau_n)H_n'(\tau_n) = \int_{\tau_n}^{x_0}\int_{\tau_n}^{y} \tilde{s}_n(u)\,du\,dy$$

and

$$Y_n(x_0) - Y_n(\tau_n) - (x_0 - \tau_n)Y_n'(\tau_n) = \int_{\tau_n}^{x_0} (U_n(u) - U_n(\tau_n))\,du.$$

Therefore, the first term can be written as

$$n^{4/5} \int_{y=\tau_n}^{x_0} \left( \int_{u=\tau_n}^{y} \tilde{s}_n(u)\,du - U_n(y) - U_n(\tau_n) \right) dy.$$

Adding and subtracting $n^{4/5} \int_{\tau_n}^{x_0} (U_0(y) - U_0(\tau_n)) \, dy = n^{4/5} \int_{\tau_n}^{x_0} \int_{\tau_n}^{y} s_0(u) \, du \, dy$, this expression can in turn be written as

$$n^{4/5} \int_{\tau_n}^{x_0} \int_{\tau_n}^{y} (\tilde{s}_n(u) - s_0(u)) \, du \, dy$$
$$- n^{4/5} \int_{\tau_n}^{x_0} (U_n(u) - U_0(u) - (U_n(\tau_n) - U_0(\tau_n))) \, du.$$

Using a second-order Taylor expansion of $s_0$ around $x_0$, this expression can be seen to equal

$$n^{4/5} \int_{\tau_n}^{x_0} \int_{\tau_n}^{y} (\tilde{s}_n(u) - s_0(x_0) - (u - x_0) s_0'(x_0)) \, du \, dy$$
$$- n^{4/5} \int_{\tau_n}^{x_0} (U_n(u) - U_0(u) - (U_n(\tau_n) - U_0(\tau_n))) \, du$$
$$- \frac{1}{2} n^{4/5} \int_{\tau_n}^{x_0} \int_{\tau_n}^{y} (u - x_0)^2 s_0''(\xi_n) \, du \, dy = O_P(1),$$

with $\xi_n \in (\tau_n, x_0)$, by Lemmas A.3, A.5 and A.6. □

LEMMA A.6. *Assume the kernel $k$ satisfies Assumption 2.6. Then*

$$Z_n(t) := n^{3/5}((U_n - U_0)(x_0 + n^{-1/5}t) - (U_n - U_0)(x_0)) \xrightarrow{d} \frac{\sqrt{g_0(x_0)}}{k(0)} W(t),$$

*in the space $D(-\infty, \infty)$ endowed with the topology of uniform convergence on compacta. Here, $W$ denotes a two-sided standard Wiener process.*

PROOF. By equation (14), we can write

$$U_n(x) = V_n(x) - \frac{1}{k(0+)} \mathbb{G}_n(x) \qquad \text{with } V_n(x) = x - \int_0^x \mathbb{G}_n(x-s) \ell(s) \, ds.$$

Define $V_0$ analogously, replacing $\mathbb{G}_n$ by $G_0$. It is easy to see that

(49) $$Z_n(t) = Z_n^{(1)}(t) - \frac{n^{3/5}}{k(0)} \int_{x_0}^{x_0 + n^{-1/5}t} d(\mathbb{G}_n - G_0)(x),$$

where

$$Z_n^{(1)}(t) = n^{3/5}((V_n - V_0)(x_0 + n^{-1/5}t) - (V_n - V_0)(x_0)).$$

The last term on the right-hand side of (49) converges to the two-sided Wiener process as indicated in the statement of the lemma. For the first term, we can write

$$Z_n^{(1)}(t) = n^{3/5} \left( \int_0^{x_0 + n^{-1/5}t} (G_0(y) - \mathbb{G}_n(y)) \ell(x_0 + n^{-1/5}t - y) \, dy \right.$$

$$- \int_0^{x_0} (G_0(y) - \mathbb{G}_n(y))\ell(x_0 - y)\, dy \bigg)$$
$$= n^{3/5} \bigg( \int_0^{x_0} (G_0(y) - \mathbb{G}_n(y))(\ell(x_0 + n^{-1/5}t - y) - \ell(x_0 - y))\, dy \bigg)$$
$$+ n^{3/5} \bigg( \int_{x_0}^{x_0 + n^{-1/5}t} (G_0(y) - \mathbb{G}_n(y))\ell(x_0 + n^{-1/5}t - y)\, dy \bigg).$$

Hence, for any $M > 0$, we get for $n$ sufficiently large that
$$\sup_{|t| \leq M} |Z_n^{(1)}(t)| \leq n^{3/5}(\|\mathbb{G}_n - G_0\|_\infty CMn^{-1/5} + 2n^{-1/5}M\|\mathbb{G}_n - G_0\|_\infty \ell(0^+))$$
$$= O_P(n^{-1/10}). \qquad \square$$

**Acknowledgements.** We thank two anonymous referees for their helpful remarks which improved the paper substantially.

Delft Institute of Applied Mathematics (DIAM)
Faculty of Electrical Engineering,
  Mathematics and Computer Science
Delft University of Technology
Mekelweg 4, 2628 CD Delft
The Netherlands
E-mail: f.h.vandermeulen@tudelft.nl